\documentclass[12pt,twoside]{article} 

\usepackage[colorlinks=false]{hyperref}
\usepackage{amssymb}
\usepackage{graphicx}
\usepackage{amscd}
\usepackage{graphics,amsmath,amssymb}
\usepackage{amsfonts}
\usepackage{latexsym}
\usepackage{epsf}

\date{} 

\begin{document} 

\newcommand{\li}{\mathop{\mathrm{li}}}
\newcommand{\median}{\mathop{\mathrm{median}}}
\newcommand{\mean}{\mathop{\mathrm{mean}}}
\newcommand{\Exp}{\mathop{\mathrm{Exp}}}
\newcommand{\Gumbel}{\mathop{\mathrm{Gumbel}}}

\centerline{\bf International Mathematical Forum, Vol.\,13, 2018, no.\,2, 65--78} 


\centerline{\bf http://dx.doi.org/10.12988/imf.2018.712103}

\centerline{} 

\centerline{} 

\centerline{} 

\centerline{\Large{\bf On the $n$th Record Gap Between Primes}} 

\centerline{} 

\centerline{\Large{\bf in an Arithmetic Progression}} 

\centerline{} 

\centerline{\bf {Alexei Kourbatov}} 

\centerline{} 

\centerline{www.JavaScripter.net/math} 

\centerline{15127 NE 24th Street \#578} 

\centerline{Redmond, WA, 98052, USA}

\newtheorem{Theorem}{\quad Theorem}[section] 

\newtheorem{Definition}[Theorem]{\quad Definition} 

\newtheorem{Corollary}[Theorem]{\quad Corollary} 

\newtheorem{Lemma}[Theorem]{\quad Lemma} 

\newtheorem{Example}[Theorem]{\quad Example} 

\centerline{}

{\footnotesize Copyright $\copyright$ 2018 Alexei Kourbatov. This article is distributed under the Creative Commons Attribution License, which permits unrestricted use, distribution, and reproduction in any medium, provided the original work is properly cited.}

\centerline{}

\begin{abstract} 
Let $q>r\ge1$ be coprime integers. Let
$R(n,q,r)$ be the $n$th record gap between primes in the arithmetic progression $r$, $r+q$, $r+2q,\ldots,$
and denote by $N_{q,r}(x)$ the number of such records observed below $x$.
For $x\to\infty$, we heuristically argue that if the limit of $N_{q,r}(x)/\log x$ exists, then the limit is 2.
We also conjecture that $R(n,q,r)=O_q(n^2)$.
Numerical evidence supports the conjectural a.s.\ upper bound
$$R(n,q,r)<\varphi(q)n^2+(n+2)q\log^2 q.$$ 
The median (over $r$) of $R(n,q,r)$ grows like a quadratic function of $n$;
so do the mean and quartile points of $R(n,q,r)$. For fixed values of $q\ge200$ and $n\approx10$,
the distribution of $R(n,q,r)$ is skewed to the right and close to both Gumbel and lognormal distributions;
however, the skewness appears to slowly decrease as $n$ increases.
The existence of a limiting distribution of $R(n,q,r)$ is an open question.
\end{abstract} 

{\bf Mathematics Subject Classification:} 11N05, 11N13, 11N56 \\

{\bf Keywords:} arithmetic progression, Cram\'er conjecture, distribution of records,
prime gap, residue class, Resnick duality theorem, Shanks conjecture.

\pagebreak
\section{Introduction} 

Let $q>r\ge1$ be coprime integers, and consider the arithmetic progression
\begin{equation}\label{eq1}
r,\ r+q,\ r+2q,\ r+3q,\ldots 
\end{equation}
Dirichlet \cite{dirichlet} proved that there are infinitely many primes in this progression. 
Let $R(n,q,r)$ be the size of the $n$th record gap between primes in the arithmetic progression (\ref{eq1}),
and denote by $N_{q,r}(x)$ the number of such records observed below $x$.
The prime number theorem for arithmetic progressions \cite{friedgold1996} guarantees that, 
for any fixed coprime pair $(q,r)$, the sequence of records $R(n,q,r)$ is infinite.
Our present goal is to investigate the behavior of $R(n,q,r)$ and $N_{q,r}(x)$ statistically,
using numerical evidence and heuristic reasoning. 

In \cite{kourbatov2016} we already studied maximal gaps between primes {\em below} $x$ in progression (\ref{eq1}).
We empirically found that, for a fixed $q$, the histogram of appropriately rescaled 
maximal gaps between primes $r+kq\le x$ is very close to the Gumbel extreme value distribution.
The nature of empirical results in \cite{kourbatov2016} is akin to probabilistic results for {\it sample maxima}.

For i.i.d.~random variables, Resnick's theorems \cite{resnick1973,resnick1973L} 
establish that the limit law of sample maxima {\it cannot be the same} as the limit law of the $n$th record.
In particular, if the limit law of sample maxima is the Gumbel distribution, 
then the {\it normal\,} distribution is a possible limit law for the $n$th record.
(We will refer to this situation as the ``Gumbel/normal case'' of Resnick's theorems.)
Three generally possible limiting distributions of the $n$th records 
in sequences of i.i.d.~random variables are the normal, lognormal, and negative 
lognormal distributions \cite{arnold1998,resnick1973L}. 
However, beyond i.i.d.\ settings, one can also encounter situations where
the Gumbel distribution itself may be the limit law for records \cite[p.\,193]{arnold1998}.

Probabilistic results \cite{arnold1998,krug2007,resnick1973,resnick1973L} are not directly applicable to record prime gaps. 
Nevertheless, it is natural to look at the growth and distribution of the record gaps 
between primes in progressions (\ref{eq1})~--- and investigate whether any results 
similar to Resnick's duality theorem might also be true for record gaps $R(n,q,r)$.
Does the number of records $N_{q,r}(x)$ behave like it would in an i.i.d.~sequence, 
i.e., does $N_{q,r}(x)$ grow about as fast as $\log x$?
What is the order of magnitude of the $n$th record gap, as a function of $n$? 
What are statistical properties of the $n$th record gap? 
Is the actual distribution of $R(n,q,r)$ approximately normal or lognormal or Gumbel~--- or none of the above?
We will attempt to answer these questions using heuristics and statistical analysis of numerical results.
Still, we have to remember that {\it prime numbers} are neither random nor independent 
\cite{granville,pintz2007}; likewise, {\it prime gaps} are neither random nor independent.
So any statistical observations and heuristic reasoning about prime gaps 
should be used with a lot of caution.

\section{Heuristic predictions}\label{heuristic}

\subsection{The $n$th record gap between primes}\label{heur-primes}

Denote by $G(x)$ the maximal gap between primes below $x$.
Let $R(n)$ be the $n$th record prime gap; $R(n)={\rm{A005250}}(n)$ in the 
{\it Online Encyclopedia of Integer Sequences} (OEIS) \cite{oeis}.
Suppose that $x$ is so large that there have been {\it many record gaps} between primes below $x$.
Cram\'er \cite{cramer} used probabilistic reasoning to conjecture that 
\begin{equation}\label{cramerconj}
G(x)=O(\log^2 x),
\end{equation}
while Shanks \cite{shanks} heuristically found 
\begin{equation}
G(x)\sim\log^2x \qquad\mbox{ as }x\to\infty \qquad\mbox{ \cite[p.\,648]{shanks}}.
\end{equation}

Let $\tau=\tau(x)$ be a function estimating the number of record prime gaps with endpoints in $[x,ex]$. 
We postulate the existence of such a function $\tau$ and, 
in accordance with our earlier observations \cite[section 3.4]{kourbatov2016}, assume that 

\medskip\noindent
(A) $\tau \ge 1$ as $x\to\infty$.
This means that prime gap records occur more often 
than records in a sequence of i.i.d.\ random variables, for which we would have\footnote{
 We expect $N\approx\li x$ gaps between primes below $x$, and 
 about $\log N$ records in an i.i.d.\ random sequence with $N$ terms; \ 
 $\li x$ denotes the logarithmic integral.
} 
$\log\li(ex)-\log\li x<1$ (for all $x\ge3$), 
while $\lim\limits_{x\to\infty} (\log\li(ex)-\log\li x)=1$.

\medskip\noindent
(B) $\tau = o(\log x)$ as $x\to\infty$.
Together with (\ref{cramerconj}), this means that 
only a zero proportion of positive integers are values of the $R(n)$ function.

\medskip\noindent
(C) $\tau$ is a continuous, non-decreasing, slowly varying function\footnote{
 Here we do not assume that $\tau$ tends to a finite limit as $x\to\infty$; but see sect.\,\ref{limit-N-log-x}.
}
of $x$. 

\medskip\noindent
Assume further that the actual number of records in $[x,ex]$ does not differ much from $\tau$.
Then there are about $\bar\tau\log x$ records below $x$, where $\bar\tau$ is the 
average value of $\tau$ on the interval $[1,x]$.
Denoting by $n$ the number of records observed up to $x$, we have
\begin{equation}
{n\over\bar\tau} \sim \log x.
\end{equation}
This, together with the Cram\'er and Shanks conjectures, implies that 
\begin{equation}
R(n) = G(x) \sim \log^2x \sim {n^2\over\bar\tau^2} \le n^2 \quad\mbox{ as }x\to\infty.
\end{equation}
Granville's correction \cite[p.\,24]{granville} to the Cram\'er and Shanks conjectures 
might imply an additional numerical constant 
in the above estimate; still we have $R(n)=O(n^2)$.

\medskip\noindent
{\bf Reality check.}
Computations of Oliveira e Silva, Herzog and Pardi \cite{toes2014} and,
more recently, Jacobsen, Nair, and others \cite{nicely}
established the actual size of the $n$th record prime gap $R(n)=G(x)$, for $x\le10^{19}$  and $n\le77$.
The actual prime gaps indeed turn out to satisfy
\begin{eqnarray}
G(x) &\lesssim& \log^2x \,\ \qquad\qquad\mbox{ for } x \le 10^{19}, \\
R(n) &\le& n^2        \qquad\qquad\qquad\mbox{ for } n \le 77, \\
R(n) &\approx& 0.25 n^2 + 0.5n   \,\quad\mbox{ for } n \le 77,  \label{quarter-n-sq}
\end{eqnarray}
which suggests that in the available data range we can
 take $1/\bar\tau^2\approx0.25$ and $\bar\tau\approx2$.

\subsection{The $n$th record gap $R(n,q,r)$}

Now consider the general case: 
gaps between primes in the arithmetic progression (\ref{eq1}).
Suppose that $x$ is so large that we have already observed {\it many record gaps} between primes $\le x$ 
in progression (\ref{eq1}).
Let $G_{q,r}(x)$ be the maximal gap between primes $r+kq\le x$, $k\in{\mathbb N}^0$.

Instead of the Cram\'er and Shanks conjectures, we will now need the following more general statements 
\cite[sections 5.2, 5.3]{kourbatov2016}:

\medskip\noindent
{\bf Generalized Cram\'er conjecture.}
Almost all maximal gaps $G_{q,r}(x)$ satisfy
\begin{equation}\label{gencramer}
G_{q,r}(x) ~<~ \varphi(q) \log^2 x
\end{equation}
for any coprime $q>r\ge1$. Here $\varphi(q)$ is {\it Euler's totient function}.

\medskip\noindent
{\bf Generalized Shanks conjecture.}
Almost all maximal gaps $G_{q,r}(x)$ satisfy 
\begin{equation}\label{genshanks}
G_{q,r}(x) ~\sim~ \varphi(q) \log^2 x \qquad\mbox{ as }x\to\infty.
\end{equation}

The heuristic reasoning then proceeds similar to the previous subsection. 
Let $\tau$ be a function estimating the number of record gaps 
between primes $p=r+kq$ with end-of-gap primes $p\in [x,ex]$.
As before, for any fixed $q$, let $\tau=\tau(q,x)$ obey the heuristic assumptions 
(A), (B), (C) of sect.\,\ref{heur-primes}.
There are about $\bar\tau\log x$ record gaps between primes $r+kq\le x$; 
denoting by $n$ the ``typical'' number of records up to $x$
we have
\begin{equation}\label{n-logx}
n 
\sim \bar\tau\log x \qquad\mbox{ as }x\to\infty.
\end{equation}
Using eqs.\,(\ref{gencramer}), (\ref{genshanks}), (\ref{n-logx}) we estimate the ``typical'' $n$th record gap:
\begin{equation}\label{rough}
R(n,q,r)~=~G_{q,r}(x)
 ~\lesssim~ \varphi(q) \log^2 x 
 ~\sim~  \varphi(q) {n^2\over\bar\tau^2}.
\end{equation}
The above is valid for {\it large} $n$ and $x$. 
To make estimate (\ref{rough}) applicable to moderate $n$,
we add a semi-empirical correction term of size $O_q(n)$ (motivated in part by heuristics of \cite{liprattshakan}):
\begin{equation}\label{finer-ineq}
R(n,q,r) ~\lesssim~ \varphi(q) {n^2\over\bar\tau^2} + (n+2)q\log^2 q.
\end{equation}
Roughly speaking, the correction term takes into account that 
in progression (\ref{eq1}) the very first prime $p=r+kq$ might occur unusually late\footnote{ 
  Cf.~\cite[section 2]{liprattshakan}; note
  the heuristic formula for $\limsup P(q)$ 
  describing the behavior of the first prime in progression (\ref{eq1}).
} 
and then subsequent primes occur less frequently than usual. 

But we do not have any precise knowledge of $\tau$.
Therefore, let us eliminate $\tau$ using our assumption (A) $\tau\ge1$; 
we thus heuristically arrive at the a.s.\ upper bound
\begin{equation}\label{final-ineq}
R(n,q,r) ~<~ \varphi(q)n^2 + (n+2)q\log^2 q.
\end{equation}
For large $x$, computations of \cite[section 3.4]{kourbatov2016} suggest that $\tau$ is strictly greater than one;
with this in mind, we expect at most finitely many exceptions to inequality (\ref{final-ineq}) for any fixed $q$.
In section \ref{numresults} we will compare this heuristic prediction with results of computations.

\subsection{The limit of $N_{q,r}(x)/\log x$}\label{limit-N-log-x}

As before, let $N_{q,r}(x)$ be the number of record gaps observed between primes $\le x$ in progression (\ref{eq1}).
Let $\tau=\tau(q,x)$ be an estimator for $N_{q,r}(ex)-N_{q,r}(x)$, the number of record gaps 
between primes $p=r+kq$, with $p\in [x,ex]$. As a function of $x$,
let $\tau$ obey the heuristic assumptions (A), (B), (C) of sect.\,\ref{heur-primes}.
Below we heuristically argue\footnote{
 A similar argument for gaps between prime $k$-tuples leads to the number of record gaps $N_k(x)\sim(k+1)\log x$;
 \ cf.~\cite[sect.\,4.2]{kourbatov2013}.
}
that if the limit $\lim_{x\to\infty}N_{q,r}(x)/\log x=\lim_{x\to\infty}\tau$ exists, then
the limit is 2. This will justify our earlier estimate for $N_{q,r}(ex)-N_{q,r}(x)$
given in \cite{kourbatov2016}; see eq.\,(\ref{tau2}).

Suppose that the following limits exist and are equal to some number $\tau_*$: 
$$
\lim\limits_{x\to\infty}{N_{q,r}(x)\over\log x} =
\lim\limits_{x\to\infty}{\mean_r N_{q,r}(x)\over\log x} =
\lim\limits_{x\to\infty}\tau(q,x) = \tau_* > 0.
$$
Let $n$ be a ``typical'' number of records up to $x$. 
For large $x$,  eq.\,(\ref{n-logx}) gives
\begin{equation}\label{n-taustar-logx}
n \sim \tau_* \log x.
\end{equation}
Define $\Delta R(n,q,r) = R(n+1,q,r)-R(n,q,r).$
By formula (\ref{rough}), for large $q$ and large $x$ we have 
\begin{eqnarray*}
\mean_r R(n,q,r)     &\sim& \varphi(q) {n^2\over\tau_*^2},           \\
\mean_r \Delta R(n,q,r) &=& \mean_r \big( R(n+1,q,r)-R(n,q,r) \big)  \\
                     &\sim& {\varphi(q)\over\tau_*^2} ((n+1)^2-n^2) ~\sim~ {2n\varphi(q)\over\tau_*^2}. 
\end{eqnarray*}
Combining this with (\ref{n-taustar-logx}) we find
\begin{equation}\label{mean-R-two-tau-phi-x}
\mean_r \Delta R(n,q,r) ~\sim~ {2\over\tau_*}\varphi(q)\log x.
\end{equation}
On the other hand, heuristically we expect that, on average, two consecutive record gaps
should differ by the ``local'' average gap between primes in progression (\ref{eq1}):
\begin{equation}\label{avg-gap}
\mean_r \Delta R(n,q,r) ~\sim~ \varphi(q)\log x \ \ \mbox{ (average gap near $x$)}.
\end{equation}
Together, eqs.\,(\ref{mean-R-two-tau-phi-x}) and (\ref{avg-gap}) imply that
$$
\tau_* = 2.
$$

\noindent
{\it Remark.} Computations \cite[sect.\,3.4]{kourbatov2016} yield the empirical estimate
\begin{equation}\label{tau2}
\tau(q,x) ~\approx~ 
\mean_r (N_{q,r}(ex)-N_{q,r}(x)) ~\approx~ 
2 - {\kappa(q)\over\log x - \delta(q)},  
\end{equation}
which agrees with the above heuristic prediction: 
$$\lim\limits_{x\to\infty}\tau=\tau_*=2.$$

\bigskip
{\footnotesize
\begin{center}Table~1. \ The 10th record gap between primes $r+kq$, \ $q=50$ \\[0.5em]
\begin{tabular}{rcrr}
\hline {\normalsize $\vphantom{1^{1^1}}$}
$r$  &  Record gap $R(10,50,r)$ & Start of gap  &  End of gap \\
[0.5ex]\hline
\vphantom{\fbox{$1^1$}}
 1&          1150 &  158551 &  159701   \\
 3&          1950 &  504953 &  506903   \\
 7&          1950 &  959207 &  961157   \\
 9&          1950 & 1229359 & 1231309   \\
11&          1150 &   56911 &   58061   \\
13&          1600 &  211663 &  213263   \\
17&          1400 &  404267 &  405667   \\
19&          1950 &  794669 &  796619   \\
21&          2300 & 6534071 & 6536371   \\
23&          1350 &  266023 &  267373   \\
27&          2100 & 1286777 & 1288877   \\
29&          1150 &  145879 &  147029   \\
31&          1150 &  289381 &  290531   \\
33&          3000 & 8314433 & 8317433   \\
37&          1950 & 1336637 & 1338587   \\
39&          1650 &  706039 &  707689   \\
41&          1650 & 1061591 & 1063241   \\
43&          1400 &  668543 &  669943   \\
47&\phantom{0}750 &   39847 &   40597   \\
49&          1150 &  241249 &  242399   \\
\hline
\end{tabular}
\end{center}
}

\pagebreak
\section{Numerical results}\label{numresults}

Using a modified version of PARI/GP code from \cite{kourbatov2016} we have computed 
the first fourteen record gaps $R(n,q,r)$ for all $q\le2000$.
Twenty or more records were computed for selected small values of $q$ (see e.\,g.~Fig.\,\ref{fig:fig1}).
Records $R(n,q,r)$ were also computed for selected larger values of $q$ up to 80000.
We used all admissible values of $r\in[1,q-1]$, $\gcd(q,r)=1$, 
to assemble a complete data set of record gaps for given $q$ and $n$.
(As an example, Table~1 gives the $R(n,q,r)$ data set for $n=10$, $q=50$.)
For each data set, we computed its largest and smallest values, 
mean, median, standard deviation, skewness, and quartile points.
This section summarizes our numerical results.

\subsection{Conjectural (a.s.) upper bound for $R(n,q,r)$}\label{conj-upper-bound}

All record gaps that we have computed turn out to satisfy the heuristic inequality (\ref{final-ineq}):
$$
R(n,q,r) ~<~ \varphi(q)n^2 + (n+2)q\log^2 q
$$
for all coprime $r<q\le2000, \,n<15$.

While we expect a finite number of exceptions at least for some values of $q$, thus far we have not seen any at all.
However, if we use a smaller correction term $(n+2)\varphi(q)\log^2 q$,
then there are a couple of exceptions, e.\,g.\ for $q=20$ and $q=23$.

\subsection{The growth trend of $R(n,q,r)$}\label{growth-trend}

\begin{figure}[h] 
  \centering
  \includegraphics[bb=2 3 591 412,width=5.3in,height=3.5in]{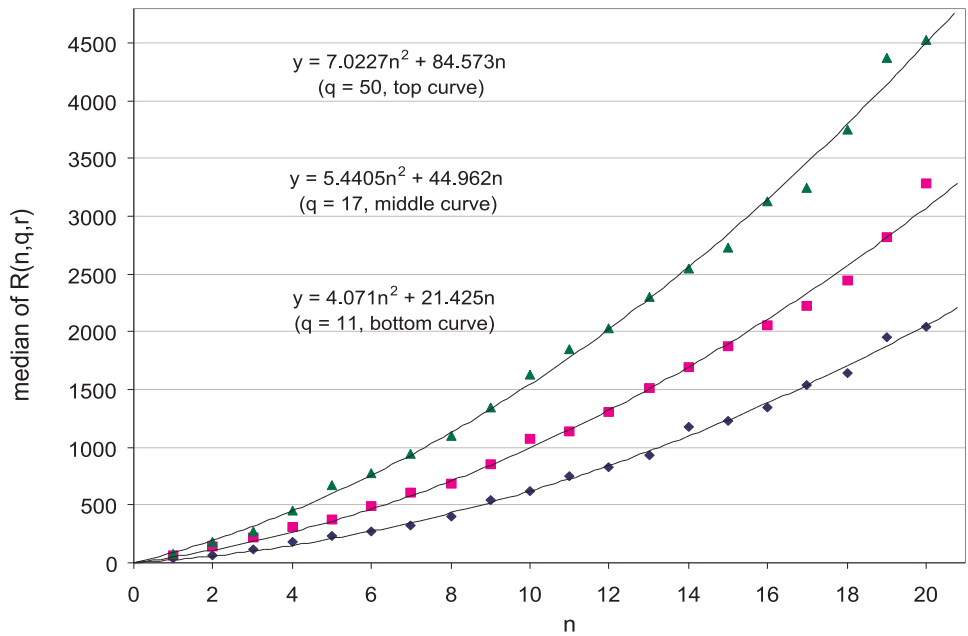}
  \caption{
    Median $n$th record gap between primes $r+kq$, \ $q=11,17,50$. \ 
    Smooth curves are quadratic approximations to median $R(n,q,r)$, eq.\,(\ref{quad-approx}).
  }
  \label{fig:fig1}
\end{figure}

For a fixed pair $(q,r)$, the sequence $R(n,q,r)$ is a strictly increasing function of $n$;
as $n$ increases, the records $R(n,q,r)$ seem to grow somewhat erratically.
But consider the {\it median} of $R(n,q,r)$ over all admissible $r$ with $\gcd(q,r)=1$. 
Table 2 and Figure \ref{fig:fig1} show that the growth of this median 
is described quite accurately by a quadratic function of $n$:
\begin{equation}\label{quad-approx}
\median\limits_{r\,\in\,[1,q]\atop q,r \,\mbox{\tiny coprime}} R(n,q,r) ~\approx~ A_q n^2 +B_q n.
\end{equation}
Rough empirical estimates for the coefficients $A_q$ and $B_q$ in (\ref{quad-approx}) are
\begin{eqnarray}
A_q &\approx& 0.3 \varphi(q), \\
B_q &<& \varphi(q)\log^2 q.
\end{eqnarray}

\pagebreak

{\footnotesize
\begin{center}Table~2. \ Median $n$th record gap between primes $r+kq$, \ $q=11,17,50$ \\[0.5em]
\begin{tabular}{rrrr}
\hline {\normalsize $\vphantom{1^{1^1}}$}
$n$  &  median $R(n,11,r)$ & median $R(n,17,r)$  &  median $R(n,50,r)$ \\
[0.5ex]\hline
\vphantom{\fbox{$1^1$}}
 1 &   33 &   68 &   75 \\
 2 &   66 &  136 &  175 \\
 3 &  110 &  221 &  275 \\
 4 &  176 &  306 &  450 \\
 5 &  231 &  374 &  675 \\
 6 &  275 &  493 &  775 \\
 7 &  319 &  612 &  950 \\
 8 &  407 &  680 & 1100 \\
 9 &  539 &  850 & 1350 \\
10 &  616 & 1071 & 1625 \\
11 &  748 & 1139 & 1850 \\
12 &  825 & 1309 & 2025 \\
13 &  935 & 1513 & 2300 \\
14 & 1177 & 1700 & 2550 \\
15 & 1232 & 1870 & 2725 \\
16 & 1342 & 2057 & 3125 \\
17 & 1540 & 2227 & 3250 \\
18 & 1639 & 2448 & 3750 \\
19 & 1958 & 2822 & 4375 \\
20 & 2046 & 3281 & 4525 \\
\hline
\end{tabular}
\end{center}
}

\begin{figure}[htb] 
  \centering
  \includegraphics[bb=10 1 576 354,width=5.8in,height=3.5in,keepaspectratio]{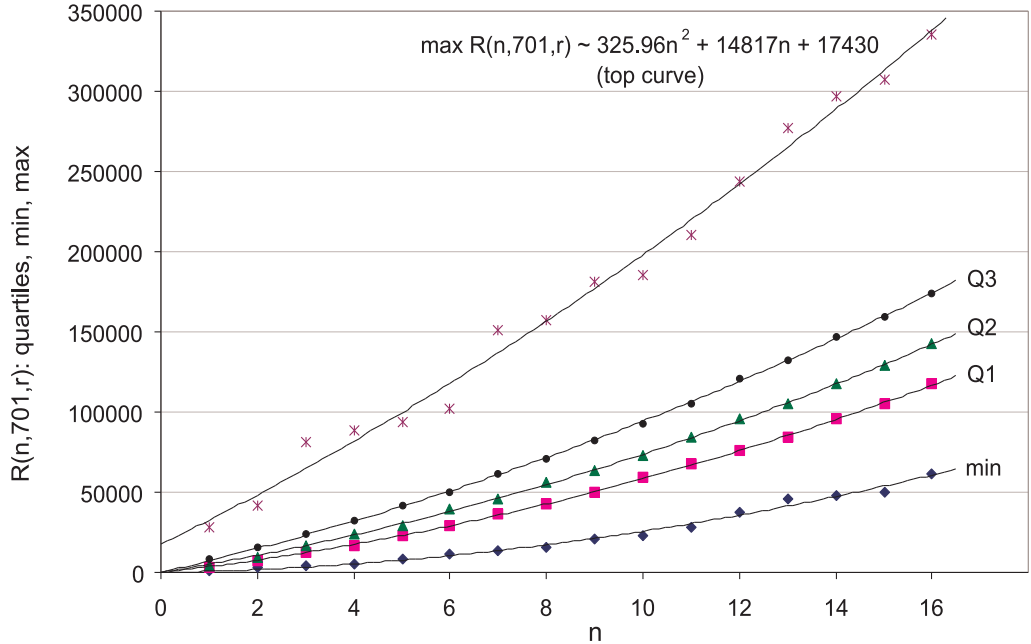}
  \caption{
    Quartile points, smallest and largest values of $R(n,q,r)$ for $q=701$. \ 
    $Q1$:~lower quartile; \ $Q2$:~median; \ $Q3$:~upper quartile. \ 
    Smooth curves are quadratic approximations.
  }
  \label{fig:fig2}
\end{figure}

Quadratic approximations similar to (\ref{quad-approx}) also work quite well 
for the mean value, least value and quartile points of $R(n,q,r)$, as shown in Figure~\ref{fig:fig2}.
For the largest value of the $n$th record, a three-term quadratic approximation is suitable
(Fig.\,\ref{fig:fig2}, top curve):
\begin{equation}\label{quad-approx-3term}
\max\limits_{r\,\in\,[1,q]\atop q,r \,\mbox{\tiny coprime}} \negthinspace R(n,q,r) 
~\approx~ \alpha_q n^2 +\beta_q n +\gamma_q.
\end{equation}

\smallskip\noindent
{\it Remark.} The quadratic approximation (\ref{quad-approx}) 
appears to remain valid for large $n$, 
with the leading coefficient $A_q$ stabilizing near some positive constant value;
this suggests that the function $\tau$ of section \ref{heuristic} tends to a finite limit as $x\to\infty$.
On the other hand, if the coefficient $A_q$ in (\ref{quad-approx}) 
were to decrease to zero when we attempt to approximate $R(n,q,r)$ for larger and larger $n$, 
this would mean that $\tau$ increases without bound. 
(In the special case of {\it record prime gaps} $R(n)$, the quadratic approximation (\ref{quarter-n-sq}) 
remains valid for a wide range of $n$, at least up to $n=77$, 
which suggests that $\lim\limits_{x\to\infty}\tau$ does exist; and the limit might be about 2.)

\subsection{The distribution of $R(n,q,r)$}\label{distsect}

\begin{figure}[hp] 
  \centering
  \includegraphics[bb=15 0 420 526,width=5.6in,height=7.1in]{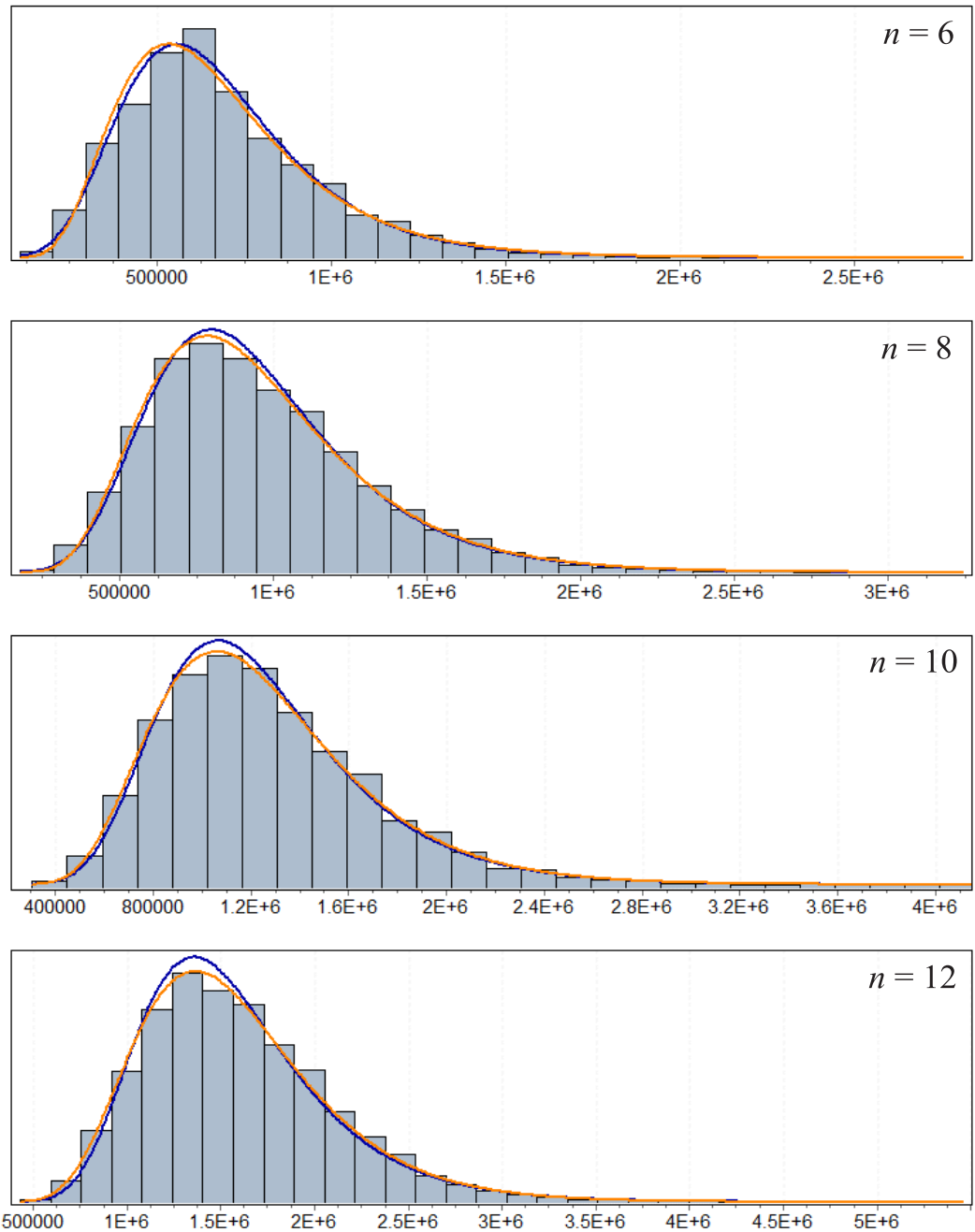}
  \caption{
    Histograms of the $n$th record gap $R(n,q,r)$ between primes in progression (\ref{eq1})
    for $q=9001$, \ $n=6,8,10,12$. \ 
    Orange curve: best-fit lognormal pdf; \ 
    dark blue curve: best-fit Gumbel pdf. \ 
  }
  \label{fig:fig3}
\end{figure}

In the previous section we have seen that, for a fixed $q$, 
the median value of $R(n,q,r)$ grows like a quadratic function of $n$.
Now let us look at the distribution of the $R(n,q,r)$ values around the median.
Figure~\ref{fig:fig3} shows the histograms of $R(n,q,r)$ computed for $q=9001$, $n=6,8,10,12$.
The histograms are clearly skewed to the right.
We see that for moderate values of $n$ both the Gumbel and lognormal distributions
are good approximations for the $R(n,q,r)$ histograms.

\begin{figure}[t] 
  \centering
  \includegraphics[bb=4 0 432 277,width=5in,height=3.2in]{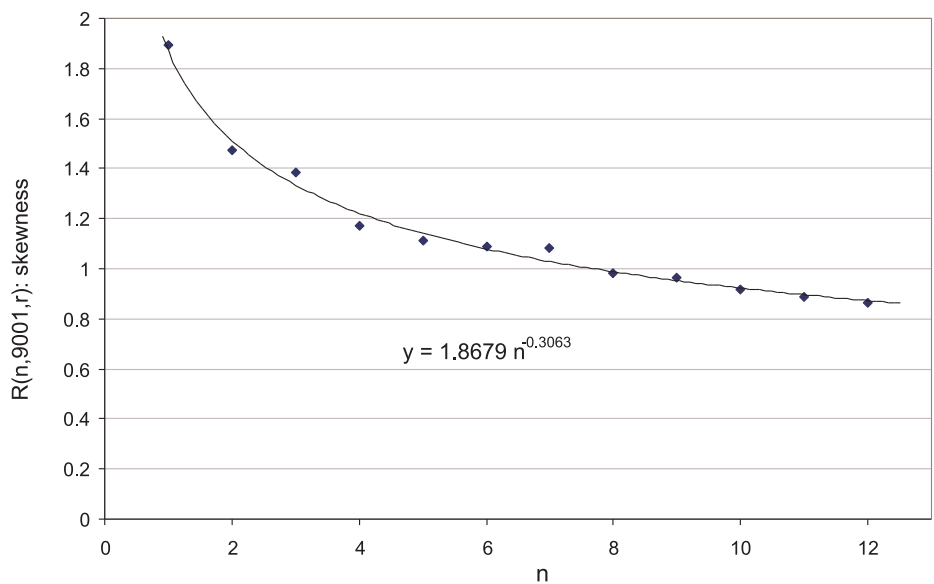}
  \caption{Skewness of the $R(n,9001,r)$ data sets for $n=1,2,3,\ldots,12$.}
  \label{fig:fig4}
\end{figure}

However, the actual $R(n,q,r)$ data sets appear to have 
slowly decreasing skewness as $n$ increases (see Fig.\,\ref{fig:fig4}),
whereas the Gumbel distribution has constant skewness independent of the distribution's scale and mode:
$$
\mbox{Gumbel distribution skewness}  ~=~ {12\sqrt{6}\,\zeta(3)\over\pi^3} ~=~ 1.139547\ldots
$$
In this respect, the lognormal distribution is a better fit to the data. 
Indeed, the best-fit lognormal distributions do reflect the decreasing skewness observed in the data.

\medskip\noindent
{\it Remark.} 
For smaller $q$, the skewness of $R(n,q,r)$ data exhibits a lot of fluctuations.
Such fluctuations may mask the general trend of decreasing skewness; 
nevertheless, this trend becomes apparent for larger $q$.

\medskip
The existence of a limiting distribution of $R(n,q,r)$ is an open question.
The decrease in skewness of $R(n,q,r)$ is satisfactorily described by a power law (Fig.\,\ref{fig:fig4}).
If the skewness continues to decrease all the way to zero, then it is possible
that the {\it normal distribution} turns out to be the limit law for $R(n,q,r)$ as $n\to\infty$.
(A sequence of lognormal distributions with vanishing skewness becomes indistinguishable from the normal distribution.)
So a certain analog of the ``Gumbel/normal case'' of Resnick's theorems \cite{resnick1973,resnick1973L} 
might be valid for the $R(n,q,r)$ limit law (if one exists);
however, convergence of records to the hypothetical normal limit law is exceedingly slow. 
For example, at the rate shown in Fig.\,\ref{fig:fig4}, 
we would need to observe over ten thousand consecutive records 
in order for the skewness to go down to 0.1.
Anyway, for moderate values of $n$ practically attainable in computation, 
the $R(n,q,r)$ histograms are quite far from the normal distribution
suggested by the ``Gumbel/normal case'' of Resnick's theorems.

The hypothetical normal limiting distribution would have a median 
that grows quadratically with $n$, as we have seen in Figures \ref{fig:fig1} and \ref{fig:fig2}.
Our computations also showed that the {\it standard deviation} of the $R(n,q,r)$ data sets
grows approximately linearly with $n$. This should be taken into account in the rescaling 
transformation to reduce the limiting distribution to standard normal form 
(if it indeed exists).

\pagebreak

\medskip
{\bf Acknowledgements.} 
I am grateful to the anonymous referees for their time and attention to this manuscript.
Thanks are also due to all contributors and editors
of the websites {\it OEIS.org} and {\it PrimePuzzles.net}.

{\bf Received: January 1, 2018; \ Published: January 22, 2018}

\end{document}